\documentclass{ps}

\usepackage{amsmath,amssymb,amsthm,amsfonts}
\usepackage[english]{babel}
\usepackage{bm}
\usepackage[numbers]{natbib}
\newcommand{\ub}{\mathbf{u}}
\newcommand{\ubhat}{\hat{\mathbf{u}}}
\newcommand{\xb}{\mathbf{x}}
\newcommand{\xbhat}{\hat{\mathbf{x}}}
\newcommand{\mb}{\mathbf{m}}
\newcommand{\mbhat}{\hat{\mathbf{m}}}
\newcommand{\fb}{\bm{f}}
\newcommand{\gb}{\mathbf{g}}
\newcommand{\xib}{\bm{\xi}}
\newcommand{\xibhat}{\hat{\bm{\xi}}}
\newcommand{\thetab}{\bm{\theta}}
\newcommand{\thetabhat}{\hat{\bm\theta}}

\newcommand{\Thetab}{\bm{\Theta}}
\newcommand{\Gb}{\mathbf{G}}
\newcommand{\Gbhat}{\mathbf{\hat{G}}}
\newcommand{\Ab}{\mathbf{A}}
\newcommand{\Abhat}{\mathbf{\hat{A}}}
\newcommand{\Bb}{\mathbf{B}}
\newcommand{\Bbhat}{\hat{\mathbf{B}}}
\newcommand{\Pb}{\mathbf{P}}
\newcommand{\Pbhat}{\hat{\mathbf{P}}}
\newcommand{\Qb}{\mathbf{Q}}
\newcommand{\Qbhat}{\hat{\mathbf{Q}}}
\newcommand{\Mb}{\mathbf{M}}
\newcommand{\Mbhat}{\hat{\mathbf{M}}}
\newcommand{\Ib}{\mathbf{I}}
\newcommand{\Fb}{\bm{F}}
\renewcommand{\hat}{\widehat}
\newcommand{\Ucal}{\mathcal{U}}
\newcommand{\Mcal}{{\mathcal M}}

\newtheorem{assumption}[thrm]{Assumption}

\begin{document}
%%-----------------------------
%%      the top matter
%%-----------------------------
\title{Consistency of direct integral estimator for partially observed systems of ordinary differential equations linear in the parameters}
\runningtitle{Consistency of direct integral estimator for partially observed ODEs}
\thanks{This research is supported by the Dutch Technology Foundation STW, which is part of the Netherlands Organisation for Scientific Research (NWO), and which is partly funded by Ministry of Economic Affairs. Part of this research was done during a visit of the second author to The Netherlands, supported by STAR Visitor Grant.}
%\thanks{...}% At most 5 thanks
%
\author{Ivan Vuja\v{c}i\'c}\address{Department of Mathematics, VU University Amsterdam, Room S3.30, De Boelelaan 1081a, 1081HV Amsterdam, The Netherlands }
\author{Itai Dattner}\address{Department of Statistics, University of Haifa, 199 Aba Khoushy Ave. Mount Carmel, Haifa 3498838, Israel}
%\author{...}\address{...}
%
%\date{...}
%
\begin{abstract} Dynamic systems are ubiquitous in nature and are used to model many processes in biology, chemistry, physics, medicine, and engineering. In particular, systems of ordinary differential equations are commonly used for the mathematical modelling of the rate of change of dynamic processes. In many practical applications, the process can only be partially measured, a fact that renders estimation of parameters of the system extremely challenging. Recently, a 'direct integral estimator' for partially observed systems of ordinary differential equations was introduced. The practical performance of the integral estimator was demonstrated, but its theoretical properties were not derived. In this paper we use the sieve framework to prove that the estimator is consistent. \end{abstract}
%
%\begin{resume} ... \end{resume}
%
\maketitle

%%-----------------------------
%%      your text
%%-----------------------------
%%%%%%%%%%%%%%%%%%%%%%%%%%%%%%%%%%%%%%%%%%%%%%%%%%%%%%%%%%%%%%%%%%%%%%%%%%%%%%%%%%%%%%%%%%%%%%%%%%%%%%%%%%%%%%%%%%%
%NEW SECTION                                    NEW SECTION                               NEW SECTION                                    NEW SECTION                   NEW SECTION                                    NEW SECTION                           NEW SECTION                                    
%%%%%%%%%%%%%%%%%%%%%%%%%%%%%%%%%%%%%%%%%%%%%%%%%%%%%%%%%%%%%%%%%%%%%%%%%%%%%%%%%%%%%%%%%%%%%%%%%%%%%%%%%%%%%%%%%%%
\section{Introduction}
\label{sec:introduction}

Mathematical models defined by a system of ordinary differential equations (ODEs) are commonly used for modelling dynamic processes. The process of interest is usually modelled by the system 
\begin{equation}\label{eq:ode_model}
\bigg\{
\begin{array}{ll}
\xb^{\prime}(t)=\Fb(\xb(t);\thetab),\ t\in [0,T],\\
\xb(0)=\xib,
\end{array}
\end{equation}
where $\xb(t)\in \xR^d $, $\xib\in\Xi\subset \xR^d$, and $\thetab\in\Thetab\subset\xR^p$. 

Given the
values of $\xib$ and $\thetab$, we denote the solution of (\ref{eq:ode_model}) by 
$$\xb(t)=\xb(t; \thetab,\xib),\ t\in [0,T].$$
The aim is to estimate the unknown parameter $\thetab$ (and if necessary $\xib$) from noisy observations
\begin{equation}
\label{observations}
Y_j(t_i)=x_j(t_i;\thetab,\xib)+\varepsilon_{j,i},\ i=1,\ldots,n, j=1,\dots,r,
\end{equation}
where $0\leq t_1<\cdots< t_n=T<\infty$ and $\varepsilon_{j,i}$ is the unobserved measurement error for $x_j$ at time $t_i$.
When the number of measured states $r$ is equal to $d$ the system is fully observed, the case well studied in the literature. Partially observed systems, i.e. when $r<d$, are common in practice but much less studied.
For an approach to estimation in this setting see \cite{ramsay2007parameter}, and \cite{qi2010asymptotic} for its asymptotic analysis. Recently, \cite{dattner2015model} extended the direct integral estimator, defined for fully observed systems linear in functions of the parameters \citep{vujavcictime,dattner2015optimal}, to partially observed ones. In the case of systems linear in the parameters, the ODE system studied in \cite{dattner2015optimal} has the form
\begin{equation}\label{eq:lin_par}
\Fb(\xb(t);\thetab)=\gb(\xb(t))\thetab,
\end{equation}
where $\gb:\xR^d\rightarrow \xR^{d\times p}$ maps the $d$-dimensional column vector $\xb$ into a $d\times p$ matrix. The introduced estimator   for fully observed  systems, i.e. when $r=d$, is motivated by the system of integral equations 
\begin{equation}\label{eq:integral_system}
\xb(t)=\xib+\int_0^t\gb(\xb(s))\xdif s \thetab,\ t\in [0,T],
\end{equation}
which follows from \eqref{eq:ode_model} and \eqref{eq:lin_par} by integration. In view of \eqref{eq:integral_system}, the estimators of the parameters $\thetab$ and $\xib$ are obtained by minimizing
$$\int_0^T\|\xbhat(t)-\xib-\int_0^t\gb(\xbhat(s))\xdif s\thetab\|^2\xdif t$$
with respect to $\thetab$ and $\xib$, where $\xbhat(t)$, $t\in[0,T]$, is a specific estimator of $\xb(t;\thetab,\xib)$. Here and subsequently, $\parallel \cdot \parallel$ denotes the Euclidean norm. Since the objective function in the display above is quadratic in $(\thetab,\xib)$, it has a unique point of minimum
 $(\thetabhat_n,\xibhat_n)$ given by
 \begin{equation}
 \label{eq:hat_parameters_full}
\begin{array}{l}
\hat{\xib}_n = \big(T\Ib_d - \hat{\Ab} \hat{\Bb}^{-1} \hat{\Ab}^{\top}\big)^{-1} \int_0^T\big\{\Ib_d -\hat{\Ab}\hat{\Bb}^{-1}\hat{\Gb}^{\top}(t)\big\} \hat{\xb}(t)\xdif t,\\
\\
\thetabhat_n=\Bbhat^{-1}\int_0^T\Gbhat^{\top}(t) \{\xbhat(t)-\xibhat_n\}\xdif t,
\end{array}
\end{equation}
where $\Ib_d$ denotes the  $d\times d$ identity matrix and
\begin{align}\nonumber
&\Gbhat(t)=\int_0^t\gb(\xbhat(s))\xdif s,\\ \label{eq:hat_terms_full}
&\Abhat =\int_0^T\hat{\Gb}(t)\xdif t,\\\nonumber
&\Bbhat=\int_0^T\Gbhat^{\top}(t)\Gbhat(t)\xdif t.\nonumber
\end{align}
If $\xbhat(\cdot)$ is a consistent estimator of $\xb(\cdot)$ in the sup norm then, under certain conditions, $(\thetabhat_n,\xibhat_n)$ is a consistent estimator of  $(\thetab,\xib)$ \citep{dattner2015optimal} .
\par 
We now describe the construction of the estimator for partially observed systems, i.e. when $r<d$, developed in \cite{dattner2015model}. Let $\Mcal$ and $\Ucal$ denote the sets of $r$-dimensional and $(d-r)$-dimensional vector functions on $[0,T]$ that correspond to $\mb(\cdot)=\mb(\cdot,\thetab,\xib)$ and $\ub^*(\cdot)=\ub^*(\cdot,\thetab,\xib)$, the measured and unmeasured components, respectively. We first construct an estimator $\mbhat_n(\cdot)$ of $\mb(\cdot)$ using the observations \eqref{observations} and for a given  $\ub\in\Ucal$,  in view of \eqref{eq:hat_terms_full}, we define
\begin{equation} \label{eq:hat_terms}
\begin{array}{l}
\xbhat_{\ub}(t)=(\mbhat_n(t),\ub(t)),\\
\\
\Gbhat_{\ub}(t)=\int_0^t\gb(\xbhat_{\ub}(s))\xdif s,\\
\\
\Abhat_{\ub} =\int_0^T\hat{\Gb}_{\ub}(t)\xdif t,\\
\\
\Bbhat_{\ub}=\int_0^T\Gbhat^{\top}_{\ub}(t)\Gbhat_{\ub}(t)\xdif t.\\
\\
\end{array}
\end{equation}
According to \eqref{eq:hat_parameters_full}, the direct integral estimator based on $\xbhat_{\ub}(\cdot)$  is
\begin{equation}
 \label{eq:hat_parameters}
\begin{array}{l}
\hat{\xib}_{\ub} = \big(T\Ib_d - \hat{\Ab}_{\ub} \hat{\Bb}_{\ub}^{-1} \hat{\Ab}_{\ub}^{\top}\big)^{-1} \int_0^T\big\{\Ib_d -\hat{\Ab}_{\ub} \hat{\Bb}_{\ub}^{-1}\hat{\Gb}^{\top}_{\ub}(t)\big\} \hat{\xb}_{\ub}(t)\xdif t,\\
\\
\thetabhat_{\ub}=\Bbhat_{\ub}^{-1}\int_0^T\Gbhat^{\top}_{\ub}(t) \{\xbhat_{\ub}(t)-\xibhat_{\ub}\}\xdif t.
\end{array}
\end{equation}
In case $\ub^*$ is measured, we can construct $\ubhat_n(\cdot)$ and the above estimator reduces to the one from \eqref{eq:hat_parameters_full}. Otherwise, for a subset $\Ucal_n\subset\Ucal$ let
\begin{equation}
\label{eq:uhat}
\ubhat_n:={\rm argmin}_{\ub\in\Ucal_n}M_n(\ub),
\end{equation}
where 
\begin{equation}
\label{eq:M_n}
M_n(\ub)=\int_0^T\|\xbhat_{\ub}(t)-\xibhat_u-\Gbhat_{\ub}(t)\thetabhat_{\ub}\|^2\xdif t.
\end{equation}
The estimators of initial value $\xib$ and the parameter $\thetab$ are 

\begin{equation}\label{eq:estimator}
\begin{array}{l}
\xibhat_n:=\xibhat_{\ubhat_n},\\
\thetabhat_n:=\thetabhat_{\ubhat_n}.
\end{array}
\end{equation}

The aim of this paper is to establish consistency of $(\thetabhat_n,\xibhat_n)$ given by (\ref{eq:estimator}) under suitable assumptions.\par
As pointed out above, \cite{dattner2015optimal} studied fully observed systems and showed that if the estimator $(\mbhat_n(\cdot),\ubhat_n(\cdot))$ of $\xb(\cdot)$ is consistent in the sup norm, then $(\thetabhat_n,\xibhat_n)$ is a consistent estimator of  $(\thetab,\xib)$. However, in partially observed systems only $\mb$ is observed. Therefore, based on the data we are only able to construct consistent estimator $\mbhat_n(\cdot)$ of $\mb(\cdot)$. Consequently, our aim is to prove that consistent $\mbhat_n(\cdot)$ gives rise to consistent estimator $\ubhat_n(\cdot)$ of $\ub^*(\cdot)$, defined in \eqref{eq:uhat}. To this end, define the deterministic counterparts of (\ref{eq:hat_terms})-(\ref{eq:hat_parameters}), 
\begin{align*}
&\xb_{\ub}(t)=(\mb(t),\ub(t)),\\
&\Gb_{\ub}(t)=\int_0^t\gb(\xb_{\ub}(s))\xdif s,\\
&\Ab_{\ub} =\int_0^T\Gb_{\ub}(t)\xdif t,\\
&\Bb_{\ub}=\int_0^T\Gb^{\top}_{\ub}(t)\Gb_{\ub}(t)\xdif t,\\
&\xib_{\ub} = \big(T\Ib_d - \Ab_{\ub} \Bb_{\ub}^{-1} \Ab_{\ub}^{\top}\big)^{-1} \int_0^T\big\{\Ib_d -\Ab_{\ub} \Bb_{\ub}^{-1}\Gb^{\top}_{\ub}(t)\big\} \xb_{\ub}(t)\xdif t,\\
&\thetab_{\ub}=\Bb_{\ub}^{-1}\int_0^T\Gb^{\top}_{\ub}(s) \{\xb_{\ub}(s)-\xib_{\ub}\}\xdif s,
\end{align*}
and the asymptotic criterion
$$M(\ub)=\int_0^T\|\xb_{\ub}(t)-\xib_{\ub}-\Gb_{\ub}(t)\thetab_{\ub}\|^2\xdif t,$$
corresponding to \eqref{eq:M_n}. Since $\ubhat_n(\cdot)$ results from a minimization over sieves $\Ucal_n$, it is a sieve extremum estimator \citep{chen2007large}. Corollary 2.6 from \cite[p. 467]{white1991some} provides the following conditions (c.f. Section 3.1, \cite[ p. 5589, 5590]{chen2007large}) which are sufficient for consistency of $\ubhat_n(\cdot)$.
\begin{itemize}
\item[{\bf C1}] $\Ucal_n\subset \Ucal_{n+1}\subset\Ucal$ and for any $\ub \in \Ucal$ there exists $\pi_n \ub\in\Ucal_n$ such that $\|\ub-\pi_n \ub\|_{\infty}\rightarrow 0$ as $n\rightarrow\infty$.
\item[{\bf C2}] $\Ucal_n$ is compact under $\|\cdot\|_{\infty}$.
\item[{\bf C3}] Functional $M(\ub)$ is continuous at $\ub^*$ in $\Ucal$  under $\|\cdot\|_{\infty}$ and $M(\ub^*) < +\infty$.
\item[{\bf C4}] $M_n(\ub)$ is a measurable function of the data $\{Y_j(t_i)\}_{j,i}$ for all $\ub\in\Ucal_n$.
\item[{\bf C5}] For any data $\{Y_j(t_i)\}_{j,i}$, $M_n(\ub)$ is lower semi-continuous on $\ub\in\Ucal_n$ under $\|\cdot\|_{\infty}$.
\item[{\bf C6}] For all $\epsilon > 0$, $M(\ub^*)< \inf_{\{\ub\in\Ucal:\|\ub-\ub^*\|_{\infty}\geq \epsilon\}} M(\ub)$.
\item[{\bf C7}]  $\sup_{\ub\in \Ucal_n}|M_n(\ub)-M(\ub)|\stackrel{P}{\rightarrow}0 $ as $n\rightarrow\infty$. %For all $n\geq 1$ 
\end{itemize}
Some additional notation: as in \cite{ding2011sieve}, for simplicity we omit the superscript $*$ in the outer probability $P^*$ whenever an outer probability applies. The metric on $\Ucal$ and $\Ucal_n $ is the one induced by $\|\cdot\|_{\infty}$. 
For a matrix function $\Mb:[0,T]\rightarrow\xR^{m\times p}$ we use the norm $\|\Mb\|_{\infty}=\sup_{t\in[0,T]}\|\Mb(t)\|$, where $\|\cdot\|$ is the Frobenius norm on $\xR^{m\times p}$. The metric on the space of matrix functions on $[0,T]$ is the one induced by $\|\cdot\|_{\infty}$. $\Mbhat$ denotes the estimator of $\Mb$ based on the data $\{Y_j(t_i)\}_{j,i}$.

 In the next section we present an example of sieves $\Ucal_n$ that satisfy conditions C1 and C2. Consistency of $\ubhat_n(\cdot)$  will be proven by verifying the conditions C3-C7. This in turn will imply consistency of $(\thetabhat_n,\xibhat_n)$. \par
The rest of the paper is organised as follows. In the next section we formulate the theorem dealing with consistency. The proofs are given in  Section 3. 
The Appendix contains technical lemmas used in Section 3.
%%%%%%%%%%%%%%%%%%%%%%%%%%%%%%%%%%%%%%%%%%%%%%%%%%%%%%%%%%%%%%%%%%%%%%%%%%%%%%%%%%%%%%%%%%%%%%%%%%%%%%%%%%%%%%%%%%%
%NEW SECTION                                    NEW SECTION                               NEW SECTION                                    NEW SECTION                   NEW SECTION                                    NEW SECTION                           NEW SECTION                                    
%%%%%%%%%%%%%%%%%%%%%%%%%%%%%%%%%%%%%%%%%%%%%%%%%%%%%%%%%%%%%%%%%%%%%%%%%%%%%%%%%%%%%%%%%%%%%%%%%%%%%%%%%%%%%%%%%%%
\section{Results}\label{previous work}
\noindent Consistency is established under the following assumptions:
\begin{assumption}
\label{assumption}
\begin{itemize}
\noindent \item[a)] (Existence and uniqueness) For any $(\thetab,\xib)\in \Theta\times\Xi$ there exists a unique solution $\xb(\cdot;\thetab,\xib)$ of \eqref{eq:ode_model} on $[0,T]$. 
\noindent \item[b)] (Identifiability) For any $(\thetab,\xib)\neq(\thetab',\xib')$ it holds that $\mb(\cdot;\thetab,\xib)\neq \mb(\cdot;\thetab',\xib')$.
\noindent \item[c)] (Densness and compactness of sieves) Sieves $\Ucal_n$ satisfy conditions {\rm C1}  and {\rm  C2}.
\noindent \item[d)] (Compactness of function spaces) $\Ucal\subset \underbrace{\xCone([0,T])\times\cdots\times \xCone([0,T])}_{d-r}$ and $\Mcal\subset \underbrace{\xCone([0,T])\times\cdots\times \xCone([0,T])}_{r}$ are compact  under $\|\cdot\|_{\infty}$ and $\mb\in\Mcal$, $\ub^*\in\Ucal$.  
\end{itemize}
\end{assumption}
Assumptions \ref{assumption} a) and b) are usual in the estimation of parameters in ODE systems, see, for example, \cite{qi2010asymptotic}. Assumption \ref{assumption} b) is equivalent to: $\mb(\cdot;\thetab,\xib)= \mb(\cdot;\thetab',\xib')\Rightarrow (\thetab,\xib)=(\thetab',\xib')$. We do not require the converse implication because it is contained in Assumption \ref{assumption} a). In other words, 
\begin{rmrk}
\label{remark:identifiability}
Assumptions \ref{assumption} a) and b) imply $\mb(\cdot;\thetab,\xib)= \mb(\cdot;\thetab',\xib')\iff (\thetab,\xib)=(\thetab',\xib')$ .
\end{rmrk}

Assumption \ref{assumption} c) is necessary for consistency of the sieve extremum estimator $\ubhat_n$; see previous section.
We now give an example of sieves that satisfy this assumption. Define
\begin{equation}
\nonumber
\Ucal_n=\{\ub\in\Ucal:u_j(t)=\sum_{k=1}^{K_{j,n}}\beta_{j,k}\phi_{j,k}(t),t\in[0,T],\beta_{j,k}\in\xR
;\sum_{j=1}^{d-r} \sum_{k=1}^{K_{j,n}}|\beta_{j,k}|\leq \Delta_n\},
\end{equation}
where $\{\phi_{j,k}\}$ is a given sequence of basis functions such that $\|\phi_{j,k}\|_{\infty}<+\infty$, $j=1,\ldots,d-r$ and $k=1,\ldots,K_{j,n}$. 
For C1 and C2 to hold we can take, for example, $\phi_{j,k}$ to be cubic splines and $\min_{1\leq j\leq d-r}(K_{j,n},\Delta_n)\rightarrow\infty$ as $n\rightarrow\infty$. Indeed, Lemma 1 from \cite{qi2010asymptotic} implies that the condition C1 is satisfied. The proof that C2 holds is the same like in \cite[p. 471]{white1991some}. Here we assume that $(\hat{\thetab}_{\ub},\hat{\xib}_{\ub})$ and $(\thetab_{\ub},\xib_{\ub})$ are well-defined on $\Ucal_n$, i.e., the matrix inverses that appear in their definition exist. Assumption \ref{assumption} d) is a technical one; it is essential for our proof. We now formulate our main result. 

%%%%%%%%%%%%%%%%%%%%%%%%%%%%%%%%%%%%%%%%%%%%%%%%%%%%%%%%%%%%%%%%%%%%%%%%%%%%%%%%%%%%%%%%%%%%%%%%%%%%%%%%%%%%%%%%%%%
%MAIN RESULT                              MAIN RESULT                         MAIN RESULT                              MAIN RESULT                               MAIN RESULT                              MAIN RESULT                     MAIN RESULT                                                 
%%%%%%%%%%%%%%%%%%%%%%%%%%%%%%%%%%%%%%%%%%%%%%%%%%%%%%%%%%%%%%%%%%%%%%%%%%%%%%%%%%%%%%%%%%%%%%%%%%%%%%%%%%%%%%%%%%%
\begin{thrm}
\label{th:main_theorem}
Let the model be defined by \eqref{eq:ode_model},\eqref{eq:lin_par},\eqref{eq:integral_system} with the map $\gb : \xR^d \rightarrow \xR^d\times\xR^p$ continuous. Fix $\xib\in\Xi$ and $\thetab\in\Theta$ and assume that $\xb(\cdot) = \xb(\cdot; \thetab, \xib)$ exists and is bounded on $[0, T]$, such that 
$$\|\xb\|_{\infty}<\infty.$$
Assume that Assumption \Rref{assumption} holds. Let $\mbhat_n(\cdot)$ be a consistent estimator of $\mb(\cdot) = \mb(\cdot; \thetab, \xib)$ in the supnorm, i.e.,
$$\|\mbhat_n-\mb\|_{\infty}\stackrel{P}{\rightarrow}0.$$
Then the estimators $\thetabhat_n$ and $\xibhat_n$  defined in \eqref{eq:estimator} are consistent, i.e.,
$$(\thetabhat_n,\xibhat_n)\stackrel{P}{\rightarrow} (\thetab, \xib)$$
holds as $n\rightarrow \infty$.
\end{thrm}
%%%%%%%%%%%%%%%%%%%%%%%%%%%%%%%%%%%%%%%%%%%%%%%%%%%%%%%%%%%%%%%%%%%%%%%%%%%%%%%%%%%%%%%%%%%%%%%%%%%%%%%%%%%%%%%%%%%%%%%%%

%%%%%%%%%%%%%%%%%%%%%%%%%%%%%%%%%%%%%%%%%%%%%%%%%%%%%%%%%%%%%%%%%%%%%%%%%%%%%%%%%%%%%%%%%%%%%%%%%%%%%%%%%%%%%%%%%%%
%NEW SECTION                                    NEW SECTION                               NEW SECTION                                    NEW SECTION                   NEW SECTION                                    NEW SECTION                           NEW SECTION                                    
%%%%%%%%%%%%%%%%%%%%%%%%%%%%%%%%%%%%%%%%%%%%%%%%%%%%%%%%%%%%%%%%%%%%%%%%%%%%%%%%%%%%%%%%%%%%%%%%%%%%%%%%%%%%%%%%%%%
\section{Proofs}\label{proofs}
Before proving the main result, we state a lemma that gives important asymptotic relationships which are used in the proof. Some of the results below implicitly use Lemma \ref{lemma_aux} of Appendix A.
%%%%%%%%%%%%%%%%%%%%%%%%%%%%%%%%%%%%%%%%%%%%%%%%%%%%%%%%%%%%%%%%%%%%%%%%%%%%%%%%%%%%%%%%%%%%%%%%%%%%%%%%%%%%%%%%%%%
%LEMMA                                    LEMMA                               LEMMA                                    LEMMA                   LEMMA                                    LEMMA                           LEMMA                             LEMMA                          LEMMA                          
%%%%%%%%%%%%%%%%%%%%%%%%%%%%%%%%%%%%%%%%%%%%%%%%%%%%%%%%%%%%%%%%%%%%%%%%%%%%%%%%%%%%%%%%%%%%%%%%%%%%%%%%%%%%%%%%%%%
\begin{lmm}
\label{main_lemma}
Let $\Ucal$ and $\Ucal_n$ satisfy Assumption \Rref{assumption} $\mathrm{c)}$, $\mathrm{d)}$. Then as $n\rightarrow\infty$
\begin{itemize}
\item [(i)] $\sup_{\ub\in\Ucal_n}\|\xbhat_{\ub}-\xb_{\ub}\|_{\infty}=o_P(1)$, $\sup_{\ub\in\Ucal_n}\|\xb_{\ub}\|_{\infty}=O(1)$, $\sup_{\ub\in\Ucal_n}\|\xbhat_{\ub}\|_{\infty}=O_P(1)$.
\item[(ii)] $\sup_{\ub\in\Ucal_n}\|\Gbhat_{\ub}-\Gb_{\ub}\|_{\infty}=o_P(1)$, $\sup_{\ub\in\Ucal_n}\|\Gb_{\ub}\|_{\infty}=O(1)$, $\sup_{\ub\in\Ucal_n}\|\Gbhat_{\ub}\|_{\infty}=O_P(1)$.
\item[(iii)] $\sup_{\ub\in\Ucal_n}\| \xibhat_{\ub}-\xib_{\ub}\|=o_P(1)$, $\sup_{\ub\in\Ucal_n}\|\xib_{\ub}\|=O(1)$, $\sup_{\ub\in\Ucal_n}\|\xibhat_{\ub}\|=O_P(1)$.
\item[(iv)] $\sup_{\ub\in\Ucal_n}\| \thetabhat_{\ub}-\thetab_{\ub}\|=o_P(1),$ $\sup_{\ub\in\Ucal_n}\|\thetab_{\ub}\|=O(1)$, $\sup_{\ub\in\Ucal_n}\|\thetabhat_{\ub}\|=O_P(1)$.
\end{itemize}
\end{lmm}
%%%%%%%%%%%%%%%%%%%%%%%%%%%%%%%%%%%%%%%%%%%%%%%%%%%%%%%%%%%%%%%%%%%%%%%%%%%%%%%%%%%%%%%%%%%%%%%%%%%%%%%%%%%%%%%%%%%%%%%%%

%%%%%%%%%%%%%%%%%%%%%%%%%%%%%%%%%%%%%%%%%%%%%%%%%%%%%%%%%%%%%%%%%%%%%%%%%%%%%%%%%%%%%%%%%%%%%%%%%%%%%%%%%%%%%%%%%%%%%%%%%
\begin{proof} In (i)-(iv) it suffices to prove only the first two statements of each, since they imply the third by using triangle inequality.

\par (i)
For any $\ub\in\Ucal_n$,  $\|\xbhat_{\ub}-\xb_{\ub}\|_{\infty}=\|\mbhat_n-\mb\|_{\infty}$ and hence the first assertion follows from consistency of $\mbhat_n(\cdot)$. The second statement follows from the compactness of $\Ucal_n$ and 
boundedness of $\mb(\cdot)$. 
%%%%%%%%%%%
\par (ii) 
Introduce $K=\sup_{\ub\in\Ucal}\|\xb_{\ub}\|_{\infty}$, which by Assumption \Rref{assumption} d) is finite. Continuity of $\gb(\cdot)$ on the closed compact ball  $B_{K+1}=\{\xb\in\xR^d:\|\xb\|\leq K+1\}$ implies its uniform continuity on $B_{K+1}$. 
Thus for $\epsilon>0$ fixed there exists $\delta>0$ such that for any $\xb_1,\xb_2\in B_{K+1}$ if $\|\xb_1-\xb_2\|<\delta$ then  $\|\gb(\xb_1)-\gb(\xb_2)\|<\epsilon/T.$
Now we show that %for any data  $\{Y_j(t_i)\}_{j,i}$  
\begin{equation}
\label{eq:implication1}
 \|\mbhat_n-\mb\|_{\infty} < \min(1,\delta) \Rightarrow  \sup_{\ub\in\Ucal_n}\|\Gbhat_{\ub}-\Gb_{\ub}\|_{\infty} \leq \epsilon.
\end{equation}
More formally, we need to prove that $\{\omega\in\Omega:\|\mbhat_n(\omega)-\mb(\omega)\|_{\infty} < \min(1,\delta)\} \subset\{\omega\in\Omega: \sup_{\ub\in\Ucal_n}\|\Gbhat_{\ub}(\omega)-\Gb_{\ub}(\omega) \|_{\infty}\leq \epsilon\}$, where $\omega$ is an outcome and $\Omega$ is the sample space. For simplicity, we suppress the explicit dependence on $\omega$ in the notation. Fix $\ub\in\Ucal_n$, $s\in[0,T]$ and assume that $\|\mbhat_n-\mb\|_{\infty} < \min(1,\delta)$. The equality  $\|\xbhat_{\ub}-\xb_{\ub}\|_{\infty}=\|\mbhat_n-\mb\|_{\infty}$ implies that $\|\xbhat_{\ub}(s)-\xb_{\ub}(s)\| < \delta$. Also, by definition of $K$ it holds that $\xb_{\ub}(s)\in B_{K+1}$, and by triangle inequality we have $\xbhat_{\ub}(s)\in B_{K+1}$.  Consequently, by uniform continuity on $B_{K+1}$ we have
$$\|\gb(\xbhat_{\ub}(s))-\gb(\xb_{\ub}(s))\|< \epsilon/T,\ s\in [0,T].$$ 
Using the derived bound we obtain that for any $\ub\in\Ucal_n$
$$\|\Gbhat_{\ub}-\Gb_{\ub}\|_{\infty}\leq \sup_{[0,T]}\int_0^t\|\gb(\xbhat_{\ub}(s))-\gb(\xb_{\ub}(s))\|\xdif s\leq \int_0^T\|\gb(\xbhat_{\ub}(s))-\gb(\xb_{\ub}(s))\|\xdif s< T\epsilon/T=\epsilon.$$
 Since this holds for any  $\ub\in\Ucal_n$ it follows that  $\sup_{\ub\in\Ucal_n}\|\Gbhat_{\ub}-\Gb_{\ub}\|_{\infty} \leq \epsilon$. Hence, \eqref{eq:implication1} is proved.
Finally, \eqref{eq:implication1} and consistency of $\mbhat_n(\cdot)$ imply
\begin{align*}
P(\sup_{\ub\in\Ucal_n}\|\Gbhat_{\ub}-\Gb_{\ub}\|_{\infty} >\epsilon)& \leq P(\|\mbhat_n-\mb\|_{\infty}  \geq  \min(1,\delta)) \rightarrow 0,
\end{align*}
as $n\rightarrow\infty$, which proves the first statement. Continuity of $\gb$ implies its boundedness on $B_{K+1}$ i.e. there exists $C>0$ such that $\|\gb(\xb)\|\leq C$ for every $\xb\in B_{K+1}$. Fix any $\ub\in\Ucal_n$ and $s\in[0,T]$.  Since $\xb_{\ub}(s)\in B_{K+1}$ we obtain
$$\|\Gb_{\ub}\|_{\infty}=\sup_{t\in[0,T]}\|\int_0^t \gb(\xb_{\ub}(s))\xdif s\|\leq \int_0^T \|\gb(\xb_{\ub}(s))\|\xdif s\leq TC,$$
which is the second claim.
\par
%%%%%%%%%%%
\par (iii) and (iv)
We first prove that 
\begin{align}
\label{eq:A_eq}
\sup_{\ub\in\Ucal_n}\|\Abhat_{\ub}-\Ab_{\ub}\|=o_P(1), \sup_{\ub\in\Ucal_n}\|\Ab_{\ub}\|=O(1),\sup_{\ub\in\Ucal_n}\|\Abhat_{\ub}\|=O_P(1),\\
\label{eq:B_eq}
\sup_{\ub\in\Ucal_n}\|\Bbhat_{\ub}-\Bb_{\ub}\|=o_P(1),\sup_{\ub\in\Ucal_n}\|\Bb_{\ub}\|=O(1), \sup_{\ub\in\Ucal_n}\|\Bbhat_{\ub}\|=O_P(1).
\end{align}
Indeed, the first assertion in \eqref{eq:A_eq} follows from Lemma \ref{lemma_1} (iii) and the result (ii) of this lemma. The second assertion is a consequence of $\sup_{\ub\in\Ucal_n}\|\Gb_{\ub}\|_{\infty}=O(1) $ and the inequality
$$\|\Ab_{\ub}\|=\| \int_0^T \Gb_{\ub}(t)\xdif t\|\leq T\|  \Gb_{\ub}\|_{\infty}.$$
By taking into account the results (ii) of this lemma and applying Lemma \ref{lemma_1} (i)-(iii) we obtain the first equality in \eqref{eq:B_eq}. The claim $\sup_{\ub\in\Ucal_n}\|\Bb_{\ub}\|=O(1) $ follows from $\sup_{\ub\in\Ucal_n}\|\Gb_{\ub}\|_{\infty}=O(1) $ and the inequality
\begin{align*}
\|\Bb_{\ub}\|&\leq \int_0^T\|\Gb^{\top}_{\ub}(t)\|\|\Gb_{\ub}(t)\|\xdif t =\int_0^T\|\Gb_{\ub}(t)\|^2\xdif t \leq T\|\Gb_{\ub}\|^2_{\infty}.
\end{align*}

Finally, repeated application of Lemma \ref{lemma_1}  and already proven results yield (iii) and (iv).
\end{proof}
%%%%%%%%%%%%%%%%%%%%%%%%%%%%%%%%%%%%%%%%%%%%%%%%%%%%%%%%%%%%%%%%%%%%%%%%%%%%%%%%%%%%%%%%%%%%%%%%%%%%%%%%%%%%%%%%%%%%%%%%%

%%%%%%%%%%%%%%%%%%%%%%%%%%%%%%%%%%%%%%%%%%%%%%%%%%%%%%%%%%%%%%%%%%%%%%%%%%%%%%%%%%%%%%%%%%%%%%%%%%%%%%%%%%%%%%%%%%%%%%%%%
\begin{proof}[Proof of Theorem \Rref{th:main_theorem}] 
By the assumption of the theorem $\mbhat_n(\cdot)$ is a consistent estimator of $\mb(\cdot)$. If  we show that $\ubhat_n(\cdot)$ is a consistent estimator of $\ub^*(\cdot)$ then 
$(\mbhat_n(\cdot),\ubhat_n(\cdot))$ is a consistent estimator of $\xb(\cdot)$. Then by Theorem 1 of 
\cite{dattner2015optimal} $(\thetabhat_n,\xibhat_n)$ is a consistent estimator of  $(\thetab,\xib)$. Consistency of $\ubhat_n(\cdot)$  is proven by verifying the conditions C3-C7 from Section \ref{sec:introduction}. 
We have divided the proof into a sequence of lemmas. That the conditions C3-C5 are satisfied is shown in Lemma \ref{lemma_4}. C6 and C7 are proven in Lemmas \ref{lemma_5} and \ref{lemma_6}, respectively. 
\end{proof}
%%%%%%%%%%%%%%%%%%%%%%%%%%%%%%%%%%%%%%%%%%%%%%%%%%%%%%%%%%%%%%%%%%%%%%%%%%%%%%%%%%%%%%%%%%%%%%%%%%%%%%%%%%%%%%%%%%%%%%%%%

%%%%%%%%%%%%%%%%%%%%%%%%%%%%%%%%%%%%%%%%%%%%%%%%%%%%%%%%%%%%%%%%%%%%%%%%%%%%%%%%%%%%%%%%%%%%%%%%%%%%%%%%%%%%%%%%%%%
%LEMMA                                    LEMMA                               LEMMA                                    LEMMA                   LEMMA                                    LEMMA                           LEMMA                             LEMMA                          LEMMA                          
%%%%%%%%%%%%%%%%%%%%%%%%%%%%%%%%%%%%%%%%%%%%%%%%%%%%%%%%%%%%%%%%%%%%%%%%%%%%%%%%%%%%%%%%%%%%%%%%%%%%%%%%%%%%%%%%%%%
\begin{lmm}
\label{lemma_4}
\begin{itemize}
\item[(i)] Functional $M(\ub)$ is continuous at $\ub^*$ in $\Ucal$ and $M(\ub^*) < +\infty$.
\item [(ii)] $M_n(\ub)$ is measurable function of the data $\{Y_j(t_i)\}_{j,i}$  for all $\ub\in\Ucal_n$.
\item[(iii)] For any data $\{Y_j(t_i)\}_{j,i}$ , $M_n(\ub)$ is lower semicontinuous on $\ub\in\Ucal_n$ under $\|\cdot\|_{\infty}$. 
\end{itemize}
\end{lmm}
%%%%%%%%%%%
\begin{proof}
(i)\par 
 We show that the mappings  $\ub\mapsto   \xb_{\ub}$, $\ub\mapsto   \Gb_{\ub}$, $\ub\mapsto   \Ab_{\ub}$, $\ub\mapsto   \Bb_{\ub}$  are continuous at $\ub^*$.
The result will then follow by repeated application of Lemma \ref{lemma_4.0}.
Fix $\epsilon>0$ and take $\delta=\epsilon$. Then  $\|\ub-\ub^*\|_{\infty}< \delta$ implies
\begin{align*}
\|\xb_{\ub}-\xb_{\ub^*}\|_{\infty}&=\|\ub-\ub^*\|_{\infty}<\delta=\epsilon,
\end{align*}
which establishes the continuity of $\ub\mapsto   \xb_{\ub}$. 
\par
To prove continuity of $\ub\mapsto   \Gb_{\ub}$ it is sufficient to prove continuity of $\xb_{\ub}\mapsto \Gb_{\ub}$ because
the composition of continuous maps is continuous. Fix $\varepsilon>0$. Let $K=\sup_{\ub\in \Ucal}\|\xb_{\ub}\|_\infty$. Under our assumption, the solution $\mb$ of the differential equation is bounded, and also any $\ub\in\Ucal$ is bounded, hence,  $\|\xb_{\ub}\|_\infty<\infty$, for any $\ub\in\Ucal$. Now we have $K<\infty$ because by (i) the mapping $\ub\mapsto  \xb_{\ub}$ is continuous and by Assumption \Rref{assumption} d) $\Ucal$ is compact. Continuity of $\gb(\cdot)$ on $\xR^d$ implies its
uniform continuity on the compact ball $B_{K}=\{\xb\in \xR^d\,|\ \|\xb\| \leq K\}$ and consequently there exists $\delta>0$ such that for all $\xb_1,\xb_2\in B_{K}$ with $ \|\xb_1-\xb_2\|<\delta$ the inequality $\|\gb(\xb_1)-\gb(\xb_2)\|<\epsilon/T$ holds.  Finally, for any $\ub\in\Ucal$ such that $\|\xb_{\ub}-\xb_{\ub^*}\|_{\infty}< \delta$ we have
$$\|\Gb_{\ub}-\Gb_{\ub^*}\|_{\infty}\leq \sup_{[0,T]}\int_0^t\|\gb(\xb_{\ub}(s))-\gb(\xb_{\ub^*}(s))\|\xdif s\leq \int_0^T\|\gb(\xb_{\ub}(s))-\gb(\xb_{\ub^*}(s))\|\xdif s< T\epsilon/T=\epsilon.$$
 Continuity of  $\xb_{\ub}\mapsto \Gb_{\ub}$ is proven. Continuity of $\ub\mapsto   \Ab_{\ub}$ and $\ub\mapsto   \Bb_{\ub}$ follows from continuity of  $\ub\mapsto   \Gb_{\ub}$ and repeated application of Lemma \ref{lemma_4.0}.

%%%%%%%%%%%

(ii)\par
Fix $\ub\in \Ucal_n$. The mapping $\{Y_j(t_i)\}_{j,i}\mapsto M_n(\ub)$ is measurable as a composition of measurable mappings. Indeed, by definition the estimator
$\mbhat_n$ is a measurable function of the data$\{Y_j(t_i)\}_{j,i}$. Also, $M_n$ is a measurable function of $\mbhat_n$ because it is continuous on $\Mcal$ under $\|\cdot\|_{\infty}$. The proof of the last claim is the same like the proof of continuity of $M$, presented in (i). 
%%%%%%%%%%%
\par
(iii)\par 
Fix any data $\{Y_j(t_i)\}_{j,i}$. Lower semicontinuity of $M_n$ is implied by its continuity. The mapping $\ub\mapsto  M_n({\ub})$ is indeed continuous  because it has the same form as 
$\ub\mapsto  M({\ub})$, with the difference that $\mb$ is substituted with $\mbhat$. But $\mbhat$ is fixed because the data is  and so the proof of continuity is the same like in (i). \end{proof}
%%%%%%%%%%%%%%%%%%%%%%%%%%%%%%%%%%%%%%%%%%%%%%%%%%%%%%%%%%%%%%%%%%%%%%%%%%%%%%%%%%%%%%%%%%%%%%%%%%%%%%%%%%%%%%%%%%%

%%%%%%%%%%%%%%%%%%%%%%%%%%%%%%%%%%%%%%%%%%%%%%%%%%%%%%%%%%%%%%%%%%%%%%%%%%%%%%%%%%%%%%%%%%%%%%%%%%%%%%%%%%%%%%%%%%%
%LEMMA                                    LEMMA                               LEMMA                                    LEMMA                   LEMMA                                    LEMMA                           LEMMA                             LEMMA                          LEMMA                          
%%%%%%%%%%%%%%%%%%%%%%%%%%%%%%%%%%%%%%%%%%%%%%%%%%%%%%%%%%%%%%%%%%%%%%%%%%%%%%%%%%%%%%%%%%%%%%%%%%%%%%%%%%%%%%%%%%%\begin{lmm}
\begin{lmm}
\label{lemma_5}
For all $\epsilon > 0$, $M(\ub^*)< \inf_{\{\ub\in\Ucal:\|\ub-\ub^*\|\geq \epsilon\}} M(\ub)$.
\end{lmm}
%%%%%%%%%%%
\begin{proof}
We will prove the statement by showing that $\ub^*$ is a unique minimum of $M$. Since for any $\ub\in\Ucal$ $M(\ub)\geq0$ and  $M(\ub^*)=0$ it follows that $\ub^*$ is a minimum of $M$. We now show that if $M(\ub)=0$ then $\ub=\ub^*$, which will imply that $\ub^*$ is the unique minimum of $M$. Fix $\ub\in\Ucal$ and assume $M(\ub)=0$. The integrand in $M$ is nonnegative and thus equal to zero Lebesgue almost everywhere. Its continuity further implies that it must be equal to zero everywhere. This yields   
\begin{align*}
\xb_{\ub}(t)-\xib_{\ub}-\int_0^t \Fb(\xb_{\ub}(s),\thetab_{\ub})\xdif s=\mathbf{0},\ t\in [0,T],\
\end{align*}
where $\mathbf{0}$ is $d$-dimensional zero vector. From the previous display we obtain that $\xb_{\ub}'(t)=\Fb(\xb_{\ub}(t),\thetab_{\ub}),\ t\in [0,T]$ and $\xb_{\ub}(0)=\xib_{\ub}$, which implies that $\xb_{\ub}(\cdot)$ is a solution of the system of ODEs
\begin{equation}
\label{eq:system_uniq}
\bigg\{
\begin{array}{ll}
\xb'(t)=\Fb(\xb(t),\thetab_{\ub}),\ t\in [0,T],\\
\xb(0)=\xib_{\ub}.
\end{array}
\end{equation}
But according to Assumption \ref{assumption} a) the solution $\xb(\cdot,\thetab_{\ub},\xib_{\ub})$ of the ODE system \eqref{eq:system_uniq} is unique, so we must have $\xb_{\ub}(\cdot)=\xb(\cdot,\thetab_{\ub},\xib_{\ub})$. This is equivalent to $\mb(\cdot,\thetab,\xib)=\mb(\cdot,\thetab_{\ub},\xib_{\ub})$ and  $\ub(\cdot)=\ub(\cdot,\thetab_{\ub},\xib_{\ub})$. By Remark \ref{remark:identifiability} it follows that 
$(\thetab,\xib)=(\thetab_{\ub},\xib_{\ub})$, which in turn, by Assumption \ref{assumption} a), implies $\xb(\cdot,\thetab,\xib)=\xb(\cdot,\thetab_{\ub},\xib_{\ub})$. Finally, from the last equality we have $\ub(\cdot,\thetab,\xib)=\ub(\cdot,\thetab_{\ub},\xib_{\ub})$, i.e. $\ub^*(\cdot)=\ub(\cdot)$. This is the desired conclusion.  
\end{proof}
%%%%%%%%%%%%%%%%%%%%%%%%%%%%%%%%%%%%%%%%%%%%%%%%%%%%%%%%%%%%%%%%%%%%%%%%%%%%%%%%%%%%%%%%%%%%%%%%%%%%%%%%%%%%%%%%%%%%%%%%%

%%%%%%%%%%%%%%%%%%%%%%%%%%%%%%%%%%%%%%%%%%%%%%%%%%%%%%%%%%%%%%%%%%%%%%%%%%%%%%%%%%%%%%%%%%%%%%%%%%%%%%%%%%%%%%%%%%%
%LEMMA                                    LEMMA                               LEMMA                                    LEMMA                   LEMMA                                    LEMMA                           LEMMA                             LEMMA                          LEMMA                          
%%%%%%%%%%%%%%%%%%%%%%%%%%%%%%%%%%%%%%%%%%%%%%%%%%%%%%%%%%%%%%%%%%%%%%%%%%%%%%%%%%%%%%%%%%%%%%%%%%%%%%%%%%%%%%%%%%%\begin{lmm}
\begin{lmm}
\label{lemma_6}
 $\sup_{\ub\in \Ucal_n}|M_n(\ub)-M(\ub)|\stackrel{P}{\rightarrow}0$, as $n\rightarrow\infty$.
\end{lmm}
\begin{proof}
We follow the idea of the proof of Proposition 3.2 of \cite{gugushvili2012n}. Inequality  $ |\|a\|^2-\|b\|^2|\leq\|a-b\|(\|a\|+\|b\|)$,  and Cauchy Shwartz and triangle inequalities in  $\xLtwo[0,T] $ imply
\begin{align*}
&|M_n(\ub)-M(\ub)|=\left|\int_0^T\left(\|\xbhat_{\ub}(t)-\xibhat_{\ub}-\Gbhat_{\ub}(t)\thetabhat_{\ub}\|^2-\|\xb_{\ub}(t)-\xib_{\ub}-\Gb_{\ub}(t)\thetab_{\ub}\|^2\right)\xdif t\right|\\
&\leq \sqrt{\int_0^T\|\xbhat_{\ub}(t)-\xb_{\ub}(t)-\xibhat_{\ub}+\xib_{\ub}-\Gbhat_{\ub}(t)\thetabhat_{\ub}+\Gb_{\ub}(t)\thetab_{\ub}\|^2\xdif t}\\
&\times \left\{\sqrt{\int_0^T\|\xbhat_{\ub}(t)-\xibhat_{\ub}-\Gbhat_{\ub}(t)\thetabhat_{\ub}\|^2\xdif t}+\sqrt{\int_0^T\|\xb_{\ub}(t)-\xib_{\ub}-\Gb_{\ub}(t)\thetab_{\ub}\|)^2\xdif t}\right\}\\
&:=\sqrt{T_1(\ub)}\{\sqrt{T_2(\ub)}+\sqrt{T_3(\ub)}\}.
\end{align*}
Since $\Ucal_n$ is compact  we have $\sup_{\ub\in\Ucal_n}T_3(\ub)=O(1)$.
Results of  Lemma \ref{main_lemma} and repeated use of Lemma \ref{lemma_1} show that $\sup_{\ub\in\Ucal_n}T_1(\ub)=o_P(1)$. Finally, from triangle inequality and inequality $(a+b)^2\leq 2(a^2+b^2)$ it follows that 
\begin{align*}
T_2(\ub)&=\int_0^T\|\xbhat_{\ub}(t)-\xibhat_{\ub}-\Gbhat_{\ub}(t)\thetabhat_{\ub}\|^2\xdif t\\
&\leq 2\int_0^T\|\xbhat_{\ub}(t)-\xb_{\ub}(t)-\xibhat_{\ub}+\xib_{\ub}-\Gbhat_{\ub}(t)\thetabhat_{\ub}+\Gb_{\ub}(t)\thetab_{\ub}\|^2\xdif t+2\int_0^T\|\xb_{\ub}(t)-\xib_{\ub}-\Gb_{\ub}(t)\thetab_{\ub}\|^2\xdif t\\
&=2T_1(\ub)+2T_3(\ub),
\end{align*}
whence $\sup_{\ub\in\Ucal_n}T_2(\ub)=O_P(1)$. This completes the proof.
\end{proof}
%%%%%%%%%%%%%%%%%%%%%%%%%%%%%%%%%%%%%%%%%%%%%%%%%%%%%%%%%%%%%%%%%%%%%%%%%%%%%%%%%%%%%%%%%%%%%%%%%%%%%%%%%%%%%%%%%%%%%%%%%

%%%%%%%%%%%%%%%%%%%%%%%%%%%%%%%%%%%%%%%%%%%%%%%%%%%%%%%%%%%%%%%%%%%%%%%%%%%%%%%%%%%%%%%%%%%%%%%%%%%%%%%%%%%%%%%%%%%
%APPENDIX                          APPENDIX                  APPENDIX              APPENDIX             APPENDIX                          APPENDIX                  APPENDIX              APPENDIX                            APPENDIX                          APPENDIX                  
%%%%%%%%%%%%%%%%%%%%%%%%%%%%%%%%%%%%%%%%%%%%%%%%%%%%%%%%%%%%%%%%%%%%%%%%%%%%%%%%%%%%%%%%%%%%%%%%%%%%%%%%%%%%%%%%%%%
\section{Appendix A}\label{sec:app}

%%%%%%%%%%%%%%%%%%%%%%%%%%%%%%%%%%%%%%%%%%%%%%%%%%%%%%%%%%%%%%%%%%%%%%%%%%%%%%%%%%%%%%%%%%%%%%%%%%%%%%%%%%%%%%%%%%%
%LEMMA                                    LEMMA                               LEMMA                                    LEMMA                   LEMMA                                    LEMMA                           LEMMA                             LEMMA                          LEMMA                          
%%%%%%%%%%%%%%%%%%%%%%%%%%%%%%%%%%%%%%%%%%%%%%%%%%%%%%%%%%%%%%%%%%%%%%%%%%%%%%%%%%%%%%%%%%%%%%%%%%%%%%%%%%%%%%%%%%%\begin{lmm}
In what follows, $\Ucal_n$ and $\Ucal$ satisfy conditions from previous sections. Recall that for a matrix function $\Mb:[0,T]\rightarrow\xR^{m\times p}$ we use the norm $\|\Mb\|_{\infty}=\sup_{t\in[0,T]}\|\Mb(t)\|$, where $\|\cdot\|$ is the Frobenius norm on $\xR^{m\times p}$. Also, $\Mbhat$ denotes the estimator of $\Mb$ based on the data $\{Y_j(t_i)\}_{j,i}$.
\begin{lmm}
\label{lemma_4.0}
For $\ub\in\Ucal_n$ $(\Ucal)$, let $\Mb_{\ub}:[0,T]\rightarrow \xR^{m\times p}$, $\Pb_{\ub}:[0,T]\rightarrow \xR^{p\times s}$  and $\Qb_{\ub}\in\xR^{m\times m}$. If the mappings
$\ub\mapsto \Mb_{\ub}$, $\ub\mapsto \Pb_{\ub}$, $\ub\mapsto \Qb_{\ub}$ are continuous on $\Ucal_n$ $(\Ucal)$ then so are $\ub\mapsto \Mb_{\ub}^{\top}$, $\ub\mapsto \Mb_{\ub}\Pb_{\ub}$, $\ub\mapsto   \int_0^T \Mb_{\ub}(t)\xdif t$, $\ub\mapsto \Qb_{\ub}^{-1}$. 
\end{lmm}
%%%%%%%%%%%%%%%%%%%%%%%%%%%%%%%%%%%%%%%%%%%%%%%%%%%%%%%%%%%%%%%%%%%%%%%%%%%%%%%%%%%%%%%%%%%%%%%%%%%%%%%%%%%%%%%%%%%%%%%%%
\begin{proof} Fix $\ub_0\in\Ucal_n (\Ucal)$. Continuity of the mappings $\ub\mapsto \Mb_{\ub}^{\top}$, $\ub\mapsto \Mb_{\ub}\Pb_{\ub}$, $\ub\mapsto   \int_0^T \Mb_{\ub}(t)\xdif t$ at $\ub_0$ follows from
\begin{align}\nonumber
&\|\Mb_{\ub}^{\top}-\Mb_{\ub_0}^{\top}\|_{\infty}=\|\Mb_{\ub}-\Mb_{\ub_0}\|_{\infty},\\
\label{eq:asymptotic}
&\|\Mb_{\ub}\Pb_{\ub}-\Mb_{\ub_0}\Pb_{\ub_0}\|_{\infty}\leq \|\Mb_{\ub}\|_{\infty}\|\Pb_{\ub}-\Pb_{\ub_0}\|_{\infty}+\|\Mb_{\ub}-\Mb_{\ub_0}\|_{\infty}\|\Pb_{\ub_0}\|_{\infty},\\
\nonumber
&\|\int_0^T \Mb_{\ub}(t)\xdif t-\int_0^T \Mb_{\ub_0}(t)\xdif t\|\leq T \|\Mb_{\ub}- \Mb_{\ub_0}\|_{\infty},
\end{align}
and continuity of the mappings $\ub\mapsto \Mb_{\ub}$ and $\ub\mapsto \Pb_{\ub}$. Continuity of $\ub\mapsto \Qb_{\ub}^{-1}$ follows from continuity of $\ub\mapsto \Qb_{\ub}$ and continuity of the matrix inversion. 
\end{proof}
%%%%%%%%%%%%%%%%%%%%%%%%%%%%%%%%%%%%%%%%%%%%%%%%%%%%%%%%%%%%%%%%%%%%%%%%%%%%%%%%%%%%%%%%%%%%%%%%%%%%%%%%%%%%%%%%%%%%%%%%%

%%%%%%%%%%%%%%%%%%%%%%%%%%%%%%%%%%%%%%%%%%%%%%%%%%%%%%%%%%%%%%%%%%%%%%%%%%%%%%%%%%%%%%%%%%%%%%%%%%%%%%%%%%%%%%%%%%%
%LEMMA                                    LEMMA                               LEMMA                                    LEMMA                   LEMMA                                    LEMMA                           LEMMA                             LEMMA                          LEMMA                          
%%%%%%%%%%%%%%%%%%%%%%%%%%%%%%%%%%%%%%%%%%%%%%%%%%%%%%%%%%%%%%%%%%%%%%%%%%%%%%%%%%%%%%%%%%%%%%%%%%%%%%%%%%%%%%%%%%%\begin{lmm}
\begin{lmm}
\label{lemma_1}
For $\ub\in\Ucal_n$ let $\Mb_{\ub}, \Mbhat_{\ub}:[0,T]\rightarrow \xR^{m\times p}$, $\Pb_{\ub},\Pbhat_{\ub}:[0,T]\rightarrow \xR^{p\times s}$ and $\Qb_{\ub}\in\xR^{m\times m}$ be such that
\begin{align*}
&\sup_{\ub\in\Ucal_n}\|\Mbhat_{\ub}-\Mb_{\ub}\|_{\infty}=o_P(1),\qquad \sup_{\ub\in\Ucal_n}\|\Pbhat_{\ub}-\Pb_{\ub}\|_{\infty}=o_P(1),\qquad \sup_{\ub\in\Ucal_n}\|\Qbhat_{\ub}-\Qb_{\ub}\|=o_P(1).
\end{align*}
\begin{itemize}
 \item [(i)] If $\sup_{\ub\in\Ucal_n}\|\Mb_{\ub}\|_{\infty}=O(1)$  then as $n\rightarrow\infty$
 $$\sup_{\ub\in\Ucal_n}\|\Mbhat_{\ub}^{\top}-\Mb_{\ub}^{\top}\|_{\infty}=o_P(1),\qquad \sup_{\ub\in\Ucal_n}\|\Mb_{\ub}^{\top}\|_{\infty}=O(1),\qquad \sup_{\ub\in\Ucal_n}\|\Mbhat_{\ub}^{\top}\|_{\infty}=O_P(1).$$ 
\item [(ii)] If $\sup_{\ub\in\Ucal_n}\|\Mb_{\ub}\|_{\infty}=O(1)$ and $\sup_{\ub\in\Ucal_n}\|\Pbhat_{\ub}\|_{\infty}=O_P(1)$ then as $n\rightarrow\infty$
$$\sup_{\ub\in\Ucal_n}\|\Mbhat_{\ub}\Pbhat_{\ub}-\Mb_{\ub}\Pb_{\ub}\|_{\infty}=o_P(1),\qquad \sup_{\ub\in\Ucal_n}\|\Mbhat_{\ub}\Pbhat_{\ub}\|_{\infty}=O_P(1).$$
\item [(iii)] $\sup_{\ub\in\Ucal_n}\| \int_0^T \Mbhat_{\ub}(t)\xdif t-\int_0^T \Mb_u(t)\xdif t\|=o_P(1),$  as $n\rightarrow\infty$.
\item [(iv)]If $\sup_{\ub\in\Ucal}\|\Qb_{\ub}\|<+\infty$ then as $n\rightarrow\infty$
 $$\sup_{\ub\in\Ucal_n}\|\Qbhat_{\ub}^{-1}-\Qb_{\ub}^{-1}\|=o_P(1),\qquad \sup_{\ub\in\Ucal_n}\|\Qb_{\ub}^{-1}\|=O(1),\qquad \sup_{\ub\in\Ucal_n}\|\Qbhat_{\ub}^{-1}\|=O_P(1).$$
\end{itemize}
\end{lmm}
%%%%%%%%%%%
\begin{proof} 
The assertions regarding the boundedness in probability in (i)-(iv)  follow from other results by using triangle inequality. As for the statements regarding boundedness we have $ \sup_{\ub\in\Ucal_n}\|\Mb_{\ub}^{\top}\|_{\infty}= \sup_{\ub\in\Ucal_n}\|\Mb_{\ub}\|_{\infty}=O(1)$. Also, by Lemma \ref{lemma_4.0} the mapping $\ub\mapsto \Qb_{\ub}^{-1}$ is continuous on $\Ucal_n$ and thus bounded. This proves $ \sup_{\ub\in\Ucal_n}\|\Qb_{\ub}^{-1}\|_{\infty}=O(1)$.\par
Now we prove the statements regarding little-O in probability. Replacing $\Mb_{\ub_0}$ and $\Pb_{\ub_0}$ by $\Mbhat_{\ub}$ and $\Pbhat_{\ub}$ in \eqref{eq:asymptotic} yields inequalities which imply (i),(ii) and (iii).
\par To prove (iv) introduce $K=\sup_{\ub\in\Ucal}\|\Qb_{\ub}\|$, which is by Assumption \Rref{assumption} d) finite. Continuity of the matrix inverse on the closed compact ball  $B_{K+1}=\{\Qb\in\xR^{m\times m}: \|\Qb\|\leq K+1\}$ implies its uniform continuity on $B_{K+1}$. 
Fix $\epsilon>0$. There exists $\delta>0$ such that for any $\Qb_1,\Qb_2\in B_{K+1}$ inequality $\|\Qb_1-\Qb_2\|< \delta $ implies  $\|\Qb_1^{-1}-\Qb_2^{-1}\| < \epsilon$.
We will show that
\begin{equation}
\label{eq:delta_ineq}
\sup_{\ub\in\Ucal_n}\|\Qb_{\ub}-\Qbhat_{\ub}\|_{\infty} <\min(1,\delta)  \Rightarrow \sup_{\ub\in\Ucal_n}\|\Qbhat_{\ub}^{-1}-\Qb_{\ub}^{-1}\| \leq \epsilon.
\end{equation}
As in proof of Lemma \ref{main_lemma} (ii), here we also suppress the explicit dependence on the outcome $\omega$ in the notation. Inequality $\sup_{\ub\in\Ucal_n}\|\Qb_{\ub}-\Qbhat_{\ub}\|_{\infty} <\min(1,\delta) $ implies that for any $\ub\in\Ucal_n$ it holds $\|\Qb_{\ub}-\Qbhat_{\ub}\|< \delta$. By definition of $K$ it holds that $\Qb_{\ub}\in B_{K+1}$, and by triangle inequality we obtain $\Qbhat_{\ub}\in B_{K+1}$. Uniform continuity now implies that $\|\Qbhat_{\ub}^{-1}-\Qb_{\ub}^{-1}\| <\epsilon$. Since this holds for any  $\ub\in\Ucal_n$ it follows that $\sup_{\ub\in\Ucal_n}\|\Qbhat_{\ub}^{-1}-\Qb_{\ub}^{-1}\| \leq \epsilon$. Therefore, \eqref{eq:delta_ineq} is proved. Finally, \eqref{eq:delta_ineq} and the assumption $\sup_{\ub\in\Ucal_n}\|\Qbhat_{\ub}-\Qb_{\ub}\|=o_P(1)$ imply
$$P(\sup_{\ub\in \Ucal_n}\|\Qbhat_{\ub}^{-1}-\Qb_{\ub}^{-1}\|>\epsilon)\leq P(\sup_{\ub\in \Ucal_n}\|\Qb_{\ub}-\Qbhat_{\ub}\|_{\infty} \geq\min(1,\delta))\rightarrow 0,$$
as $n\rightarrow\infty$, which is the first assertion. 
\end{proof} 
%%%%%%%%%%%%%%%%%%%%%%%%%%%%%%%%%%%%%%%%%%%%%%%%%%%%%%%%%%%%%%%%%%%%%%%%%%%%%%%%%%%%%%%%%%%%%%%%%%%%%%%%%%%%%%%%%%%%%%%%%
\begin{lmm}
\label{lemma_aux}

\begin{itemize}
\item[(i)] Let $\fb:[0,T]\rightarrow \xR^m$ be a vector-valued function in $\xLone[0,T]$. Then
$$\Big\|\int_0^T \fb(t)\xdif t  \Big\|\leq \int_0^T\|\fb(t)\|\xdif t.$$
\item[(ii)]
Let $\Mb\in \xR^{m\times p}$ and $\Pb\in \xR^{p\times s}$. Then 
$$\|\Mb\Pb\|\leq \|\Mb\|\|\Pb\|.$$
\end{itemize}
\end{lmm}
\begin{proof}
For (i) see  page 540 of \cite{jones2001lebesgue} and for (ii) page 550 of \cite{bernstein2009matrix}.
\end{proof} 

\begin{acknowledgement}
The first author is thankful to Bartek Knapik, Shota Gugushvili and Eduard Belitser for useful discussions.
\end{acknowledgement}
%%-----------------------------
%%      your bibliography
\bibliographystyle{plainnat}
\bibliography{references}
%%-----------------------------

\end{document}